
\input amssym.def
\input amssym

\magnification = \magstephalf
\parskip= 7.0pt

\def\R{{\Bbb R}}
\def\Q{{\Bbb Q}}

\def \no{\noindent}

\def\half{{1\over 2}}
\def\thalf{{\textstyle{\half}}}

\def\intl{\int\limits}
\def\({\left(}
\def\){\right)}
\def\lv{\left\vert}
\def\rv{\right\vert}

\def\ldz{{\zeta^\prime\over\zeta}}
\def\izT{\intl_0^T}
\def\Re{\hbox{Re}}
\def\Im{\hbox{Im}}

\centerline{\bf MEAN VALUES OF $\ldz $ AND THE GUE HYPOTHESIS}

\centerline{\phantom{o}}

\centerline{David W.~Farmer\footnote{*}{Research supported
by an NSF  Postdoctoral Fellowship and NSF grant \#DMS-9022140}}

\vskip .2in

\no {\bf 1. Introduction}

The Riemann Hypothesis states that the nontrivial zeros
of the Riemann $\zeta$-function are of the form
$\rho_j=\half + i\gamma_j$ with $\gamma_j\in {\R}$, but it
makes no assertion about the distribution of the numbers $\gamma_j$.
To study this distribution, the zero-counting function 
$N(T)=\#\{0<\gamma_j<T\,|\,\zeta(\half+\gamma_j)=0\}$ is a natural 
object to consider.  Von Mangoldt
proved $N(T)={1\over 2 \pi}T\log (T/2\pi e) 
+ O(\log T)$, and this contains much information about the 
individual $\gamma_j$.  For example,
$\gamma_j= 2\pi j/\log j + O(\log j)$, and the numbers
$\tilde{\gamma}={1\over 2\pi}\gamma\log\gamma$ have mean spacing 1.
Montgomery [M] began a finer study of the distribution of the
$\gamma_j$ by considering the pair-correlation function
$$
F(\alpha,T)=N(T)^{-1}\sum_{0<\gamma,\gamma^\prime<T} 
T^{i\alpha(\gamma-\gamma^\prime)} \,w(\gamma-\gamma^\prime) ,
$$
where $w(u)=4/(4+u^2)$.  Montgomery [M] showed 
$F(\alpha,T)= (1+o(1))T^{-2\alpha}\log T + \alpha + o(1)$, uniformly
for $0\le\alpha\le 1$, and he conjectured 
$F(\alpha,T)=1+ o(1)$, uniformly for $1\le\alpha\le  A$ for any
fixed $A>1$.  This conjecture has consequences for the
distribution of primes [G1] [GG] [GM], mean values of $S(t)$ [G2],
and mean values of the $\zeta$-function [CG] [GG].  As an example,
Goldston and Gonek [GG] show, assuming RH, if
$a\approx 1/\log T$, then
$$
\intl_0^T \lv \ldz (\thalf+a+it)\rv^2dt
\sim T\,{1-T^{-2a}\over 4a^2}+ \log^2T\intl_1^\infty
	(F(\alpha,T)-1) \,T^{-2a\alpha}d\alpha .
\eqno(1.0)
$$
Note that the first term dominates if and only if
Montgomery's conjecture is true.

Perhaps the most profound aspect of Montgomery's conjecture is that 
it says the distribution of gaps between zeros of the 
$\zeta$-function is the same as the Gaussian Unitary Ensemble 
distribution from random matrix theory.  In other words, the
gaps between zeros of $\zeta(s)$ are distributed
like the gaps between eigenvalues of large random Hermitian
matrices.  Numerical computations of Odlyzko [Odl] have found
amazingly close agreement between the actual distribution of
zeros and the distribution predicted by the GUE model.

Recently, Rudnick and Sarnak [RS] have considered higher 
correlation functions of zeros of the $\zeta$-function.  Their
result, which we describe later, strongly supports the conjecture
that all correlation functions of zeros of the $\zeta$-function
agree with those of the GUE model.  We refer to this as
the ``GUE Hypothesis,''  and we write GUE$N$ for the assumption
that the $n$-correlations of zeros of $\zeta(s)$ agree with
GUE for $2\le n\le N$.  For example, GUE2 is Montgomery's
conjecture.

In the spirit of (1.0), we prove,

\proclaim Theorem 1.  Assume RH and GUE3.  If $\,a\approx 1/\log T$ then
$$\eqalign{
\intl_{-T}^T \lv \ldz (\thalf+a+it)\rv^2 \ldz (\thalf+a+it)dt
\,&\sim T\log T \,{T^{-2a}\over 2 a^2} .\cr
}
$$

The paper is organized as follows.  In Section 2 we describe
the GUE Hypothesis and the result of Rudnick and Sarnak [RS].  
In Section~3 we describe the connection between GUE and other
conjectures, and we discuss the situation for $L$-functions.
In Section~4 we prove Theorem~1.  In Section~5 we indicate 
other integrals which can be evaluated by these techniques
and we discuss the limitations of these methods.

\noindent{\bf 2. The GUE Hypothesis}

Following Rudnick and Sarnak [RS], on RH and GUE we can evaluate
expressions of the form
$$
S_n(T,f)=
\mathop{{\sum}'}_{0<\gamma_1,\ldots,\gamma_n<T} 
	f\({L\over 2\pi}\gamma_1,\ldots,{L\over 2\pi}\gamma_n\)
$$
where ${\sum}^\prime$ means summation over distinct indices, 
$L=\log T$, and
$f$ satisfies
\itemitem{TF1.} $f(x_1,\ldots,x_n)$ is symmetric.

\itemitem{TF2.} $f( x+ (t,\ldots,t))=f( x)$  for $t\in {\R}$.

\itemitem{TF3.} $f( x)\to 0$ rapidly as $|x|\to \infty$ in
the hyperplane $\sum x_j =0$.

Condition TF1 is not really necessary because 
$S_n(T,f)=S_n(T,f_0)$,
where $f_0( x)= {1\over n!}\sum f(\sigma x)$,
where $\sigma$ runs over all permutations on $n$ letters.
Condition TF2 merely states that $f$ is a function of the differences
$x_i-x_j$.  Condition TF3 assures that various integrals will
converge.  The factor $L/2\pi$ is introduced because the 
$\tilde{\gamma}=L\gamma/2\pi$ have mean spacing~1. 

The GUE Hypothesis states that, under the conditions described above,
$$\eqalignno{
S_n(T,f)\sim \,&{1\over 2\pi}T L \intl_{{\R}^n} f( x)W_n(x)
		\delta(\overline{ x}) dx
&(2.0)
}$$
where
$$
W_n( x)=\det(K(x_i-x_j)), 
\ \ \ \ \ \ \ \ \ \ \ \ \ \ \ \ \ 
K(t)={\sin(\pi t)\over \pi t},
\ \ \ \ \ \ \ \ \ \ \ \ \ \ \ \ \
\overline{x}={x_1+\cdots+x_n\over n},
$$
and $\delta$ is the Dirac $\delta$-function.  Rudnick and Sarnak [RS] 
prove $(2.0)$ under the assumption that the Fourier transform
$\hat f(\xi)$ is supported in $\sum_j|\xi_j| < 2$.  They actually
prove the corresponding result for essentially all interesting 
$GL_r$ $L$-functions.  The only modification needed is to put $L=r\log T$
and $\hat f(\xi)$ must be supported in $\sum_j|\xi_j| < 2/r$.

\no {\bf 3. GUE, other conjectures, and L-functions}

We describe a conjecture related to GUE, and we speculate about
the connection between the $\zeta$-function and other $L$-functions.

The following conjecture is (7.4) in [F1].

\proclaim Conjecture.  Suppose $a,b,u,v$ are complex numbers
of modulus $\ll 1/\log T$, each 
having real part contained in $[A/\log T,B/\log T]$.
Then
$$
\intl_0^T {\zeta(\half + u + it)\,\zeta(\half + v - it)\over
	\zeta(\half + a + it)\,\zeta(\half + b  - it)}\,dt
	=\,T\!\(1+{\(1-T^{-(u+v)}\)} {(u-a)(v-b)\over (u+v)(a+b)}\)
	+o(T) ,
\eqno(3.0)
$$
uniformly for fixed positive $A,B$.

\break

We may differentiate both sides  with respect to any of
$a,b,u,v$ by using Cauchy's theorem and 
integrating around circles of size $1/\log T$.  This will increase
the error term by a factor of $\log T$, but the main term increases
by the same amount, so we get another meaningful formula.  For
example, differentiating
with respect to $u$ and $v$, and then setting $u=a$ and $v=b$ gives 
$$
\intl_0^T  \ldz (\thalf+a+it) \ldz (\thalf+b-it)\,dt
\sim T\,{1-T^{-(a+b)}\over (a+b)^2} .
$$
Put $b=a$ to get the left side of (1.0).  By the comment 
following (1.0), this implies GUE2.  The integral 
above is not more difficult to evaluate than the integral in (1.0),
and the methods of [GG] are easily modified to give the
above formula, on GUE2.

If we differentiate (3.0) with respect to $u,v$ and $a$, and then set
$u=a$ and $v=b$, we get
$$\eqalignno{
\intl_{0}^T \lv \ldz (\thalf+a+it)\rv^2 \ldz (\thalf+b+it)dt
\,&\sim T \log T\,{T^{-(a+b)} \over  (a+b)^2} ,
&(3.1)
}
$$
which implies Theorem~1.  This gives an interesting
consistency check between (3.0) and the GUE Hypothesis.
If the methods in Section~4 could be modified to express the
left side of (3.1) in terms of the triple-correlation measure
$W_3(x,y,z)$, then it might be possible to show that (3.1) implies
GUE3. 
Other formulas derived from (3.0) are given in (7.5) - (7.8) of~[F1].

Let $L(s)$ be the $L$-function associated to a cuspidal automorphic
representation of $GL_r/{\Q}$.  In [RS] it is shown that one
expects all of the correlation functions of the zeros of $L(s)$ to
be given by the GUE model.  They refer to this as 
the `universality of GUE.'  We now describe a general method for using
known results about $\zeta(s)$ to conjecture results about $L(s)$.
This can be seen as transferring GUE for $\zeta(s)$ to GUE for $L(s)$.

\proclaim The Scaling Principle.  For $|s-\half|\ll 1/\log t$ we have
the correspondence:
$$\eqalign{
L(s-\thalf)
\ \ \ \ \ \ \ \ \ \ \ \ \ 
&\longleftrightarrow
\ \ \ \ \ \ \ \ \ \ \ \ \ 
\zeta(r(s-\thalf))\phantom{\sum}
\cr
\cr
\sum_{n\le X^r} {a(n)\over n^s}
\ \ \ \ \ \ \ \ \ \ \ \ \ 
&\longleftrightarrow
\ \ \ \ \ \ \ \ \ \ \ \ \ 
\sum_{n\le X} {A(n)\over n^s} .
}$$

Some simple cases of this principle occur in the approximation of 
Dirichlet series by 
Dirichlet polynomials.  For example, the approximate functional
equation for $\zeta(s)$ consists of two Dirichlet polynomials
of length $T^\half$, while for a $GL_2$ $L$-function the polynomials
have length $T$.  Another example is Bombieri and Friedlander's [BF]
result that a $GL_r$ $L$-function cannot be well approximated by
a single  polynomial of length less than $T^r$.

A more interesting example can be found in mean-values near the
$\half$-line of the Dirichlet series multiplied by a 
Dirichlet polynomial.  In [F2], Theorem~1 and directly following,
two examples are given, and a version of the scaling principle
is described at the end of that paper's introductory section.

Applying the scaling principle to (3.0) we get, after
changing variables:
$$
\intl_0^T {L(\half + u + it)L(\half + v - it)\over
	L(\half + a + it)L(\half + b  - it)}\,dt
	=\,T\!\(1+{\(1-T^{-r(u+v)}\)} {(u-a)(v-b)\over (u+v)(a+b)}\)
			+o(T)
$$
One can use this formula to obtain the analog of (3.1), and then
apply the methods of [GG] to obtain the analog of (1.0), and so 
obtain GUE2 for $L(s)$.  In other words, (3.0) plus
the scaling principle implies that the pair correlation of zeros of
any $L$-function follow the GUE model.  One could also argue that
the scaling principle directly implies that the zeros of $L(s)$
behave like the (scaled) zeros of $\zeta(s)$, so GUE for $\zeta(s)$
transfers to GUE for $L(s)$.

\no {\bf 4.  Proof of Theorem~1}

The proof closely follows the first part of [GG].

From the elementary identities $|z|^2z={4\over 3}x^3 
+{4\over 3} iy^3 - {1\over 3}\overline{z}^3$
and $\overline{\zeta(s)} = \zeta(\overline{s})$ we obtain
$$\eqalignno{
\intl_{-T}^T \lv \ldz (\thalf+a+it)\rv^2 \ldz (\thalf+a+it)dt
=&{4\over 3} \intl_{-T}^T \({\Re\ldz(\thalf+a+it)}\right)^3 dt \cr
&-{1\over 3} \intl_{-T}^T \ldz(\thalf+a+it)^3 dt
+ {4i\over 3}   \intl_{-T}^T \({\Im\ldz(\thalf+a+it)}\right)^3 dt \cr
=& {8\over 3} \izT  \({\Re\ldz(\thalf+a+it)}\right)^3 dt 
-{2\over 3} \,\Re\izT  \ldz(\thalf+a+it)^3 dt.
\cr
&&(4.0)
}$$
Only the first term on the last line makes a contribution, as the
next Lemma shows.

\proclaim Lemma 2.  On RH, if $\,\half<\sigma<{3\over 4}$, then
$$
\izT \ldz(\sigma+it) dt \ll \log T \log(\sigma-\thalf),
$$
and if $\,n\ge2$,
$$
\izT \ldz(\sigma+it)^n  dt \ll {\log^n T\over (\sigma - \half)^{n-1}}.
$$

{\sl Proof}:  From the partial fraction formula [T],
$$
\ldz(s)= A - {1\over {s-1}} - 
	\half {\Gamma^\prime\over \Gamma} \({s\over 2}+1\) 
	+ \sum_\rho \({1\over s-\rho}+ {1\over \rho}\),
$$
where $A=\half\log \pi - \sum_\rho \Re\, \rho^{-1}$,  
follows easily the estimate
$$
\ldz (s) \ll {\log t\over s-\half} .
$$
Moving the path of integration to the 2-line, using the above estimate
on the horizontal segments, and integrating term-by-term on the 
new path, gives the estimates in Lemma~2.

Now we evaluate the first term in (4.0).
From the partial fraction formula and 
$$
{\Gamma^\prime\over\Gamma}(s)=\log s + O\({1\over 1+|t|}\)
$$
we get
$$
\Re\(\ldz(\thalf+a+it)\)=-\half\log{t\over2 \pi} + a\Sigma_a(t)
		+O\!\({1\over |t|+1}\) ,
\eqno (4.1)
$$
where
$$
\Sigma_a(t)=\sum_{\gamma} {1\over a^2+(t-\gamma)^2} .
$$
So,
$$
\eqalignno{
\intl_0^T \({\Re\ldz(\thalf+a+it)}\right)^3 dt =&
-{1\over 8} \izT \log^3{t\over 2\pi}\,dt 
+ {3a\over 4}\izT\Sigma_a(t)\log^2{t\over 2\pi}\,dt 
-{3a^2\over 2}\izT  \Sigma_a(t)^2\log{t\over 2\pi} \,dt \cr
&+a^3\izT \Sigma_a(t)^3\,dt 
+O\!\(\izT {\(\log^2t + a^2\Sigma_a(t)^2)\)}(t+1)^{-1}\,dt\) .&(4.2)\cr
}$$
It remains to evaluate the moments of $\Sigma_a(t)$.

\proclaim Proposition 3. Assume RH and GUE3.  
If $\,a\approx 1/\log T$ then
$$\eqalignno{
\Sigma_a(t) \ll& {\log(|t|+2)\over a^2} \cr
\intl_0^T \Sigma_a(t)\,dt =& {1\over 2a}\(T\log{T\over 2\pi}-T\)
				+O\!\({\log^2 T \log\log T }\)\cr
\intl_0^T \Sigma_a(t)^2\,dt 
\sim& {1\over 4a^2}T\log^2T\(1+{1-T^{-2a}\over 2a^2\log^2T}\) 
\cr
\intl_0^T \Sigma_a(t)^3\,dt 
\sim& {1\over 8a^3}T\log^3T\(1+ {3\over 2a^2\log^2T}\) .
}$$

By Proposition~3, and integration-by-parts, we can evaluate 
the terms in (4.2):
$$
\eqalign{
\intl_0^T \({\Re\(\ldz(\thalf+a+it)\)}\right)^3 dt \sim&
-{1\over 8}\,T\log^3T + {3\over 8}\,T\log^3T
-{3\over 8}\, T\log^3T\(1+{1-T^{-2a}\over 2a^2\log^2T}\)  \cr
&+{1\over 8}\,T\log^3T\(1+ {3\over 2a^2\log^2T}\)
\cr
\sim&\,{3\over 8}\,T\log T\,{T^{-2a}\over 2a^2} .
}$$
Combining this with (4.0) and Lemma~2 gives Theorem~1.
It remains to prove Proposition~3.  This involves putting the expressions
in Proposition~3 in a form where we can apply GUE.

{\sl Proof of Proposition~3}.
The first estimate follows from
$$
\sum_{|\gamma-t|<1}1\ll \log t,
\eqno (4.3)
$$
which follows from the zero-counting result 
$N(T)={1\over 2\pi} T\log (T/2\pi e) + O(\log T)$.

To prove the second result, combine Lemma~2 and equation (4.1):
$$\eqalign{
\intl_0^T \Sigma_a(t)\,dt  =& {1\over 2a}\intl_0^T \log{t\over 2\pi}\,dt
				+O\!\({\log T \log a\over a}\)\cr
			=&  {1\over 2a}\( T\log {T\over 2\pi} - T\)
				+O\!\({\log^2 T \log \log T}\) .
}$$

The proof of the third result is given in [GG]; this requires GUE2.
We will use essentially
the same method to evaluate the last expression in Proposition~3;
this requires GUE3.

\proclaim Lemma 4.  If $\,a\gg 1/\log T$ then
$$
\intl_0^T \Sigma_a(t)^n \, dt  = \intl_{-\infty}^\infty \(
	\sum_{0<\gamma<T} {1\over a^2+(t-\gamma)^2}
						\right)^n dt
+O\!\(\log^{3n}T \)
$$

{\sl Proof}.  The method used here is known as ``Montgomery's trick.''
We truncate the sum so that it is over $0<\gamma<T$, and
then extend the integral to go over all $t$. 
First note that
if $1<t<T$ then
$$
\sum_{\gamma\not\in [0,T]} {1\over a^2 + (t-\gamma)^2}
\ll \({1\over T+1}+{1\over T-t+1}\) {\log T\over a^2} .
$$
This follows from (4.3) and $N(T)\ll T\log T$.  
Combining this with the first estimate in Proposition~3 shows
that we may restrict all the summations to $0<\gamma<T$
with an error of $O(\log^{3n} T)$.  A similar argument shows that
we may extend the range of integration to go over $(-\infty,\infty)$
at the cost of an even smaller error.  This proves Lemma~4.

Expanding the integrand in Lemma~4 and exchanging summation and
integration gives
$$
\intl_0^T \Sigma_a(t)^3 \, dt  
\sim \sum_{0<\gamma,\gamma^\prime,\gamma^{\prime\prime}<T}
\ \intl_{-\infty}^\infty {1\over a^2 +(t-\gamma)^2}
	{1\over a^2 +(t-\gamma^\prime)^2}
		{1\over a^2 +(t-\gamma^{\prime\prime})^2} \,dt .
$$
This integral can be evaluated explicitly.
A straightforward calculation yields:

\proclaim Lemma 5.  If $\,a>0$, 
$$\eqalign{
\intl_{-\infty}^\infty &{1\over a^2 +(t-\gamma)^2}
	{1\over a^2 +(t-\gamma^\prime)^2}
	{1\over a^2 +(t-\gamma^{\prime\prime})^2} \,dt\cr
&\phantom{XXXXXXXXXXX}=
{12\pi a}F_1(\gamma,\gamma^\prime,\gamma^{\prime\prime}) + 
{\pi\over a}\(F_2(\gamma,\gamma^\prime,\gamma^{\prime\prime})+
F_2(\gamma^\prime,\gamma^{\prime\prime},\gamma)+
F_2(\gamma^{\prime\prime},\gamma,\gamma^\prime)\),
}$$
where
$$\eqalign{
F_1(x,y,z)&={1\over 4a^2+\({x-y}\right)^2}
		{1\over 4a^2+\({y-z}\right)^2}
		{1\over 4a^2+\({z-x}\right)^2}\cr
F_2(x,y,z)&={1\over 4a^2+\({x-y}\right)^2}
		{1\over 4a^2+\({y-z}\right)^2} .
}$$

This brings us to
$$\eqalignno{
\izT\Sigma_a(t)^3\,dt \sim {12\pi a}
	\sum_{0<\gamma,\gamma^\prime,\gamma^{\prime\prime}<T}
		F_1(\gamma,\gamma^\prime,\gamma^{\prime\prime})
\,+&\,
{3\pi \over a}
	\sum_{0<\gamma,\gamma^\prime,\gamma^{\prime\prime}<T}
			F_2(\gamma,\gamma^\prime,\gamma^{\prime\prime}).
&(4.4)
}$$
We will use the GUE Hypothesis to evaluate these sums.
Everything must be expressed in terms of sums over distinct indices.  
For the sum on $F_1$ we have,
$$\eqalign{
\sum_{0<\gamma,\gamma^\prime,\gamma^{\prime\prime}<T}
	F_1(\gamma,\gamma^\prime,\gamma^{\prime\prime})  =&
\sum_{0<\gamma,\gamma^\prime,\gamma^{\prime\prime}<T}
f_1\(\tilde{\gamma},\tilde{\gamma}^\prime,\tilde{\gamma}^{\prime\prime}\) \cr
=&
\mathop{{\sum}'}_{0<\gamma,\gamma^\prime,\gamma^{\prime\prime}<T}
f_1\(\tilde{\gamma},\tilde{\gamma}^\prime,\tilde{\gamma}^{\prime\prime}\)
+ {3\over 4a^2} \mathop{{\sum}'}_{0<\gamma,\gamma^\prime<T}
f_1\(\tilde{\gamma},\tilde{\gamma}^\prime\)
+{1\over 2^6a^6} \sum_{0<\gamma<T} 1
}$$
where
$$\eqalign{
f_1(x,y,z)=&{L^6\over 2^6\pi^6}
		\({L^2a^2\over \pi^2}+ (x-y)^2\)^{-1}
		 \({L^2a^2\over \pi^2}+ (y-z)^2\)^{-1}
		\({L^2a^2\over \pi^2}+ (z-x)^2\)^{-1} \cr
f_1(x,y)=&{L^4\over 2^4\pi^4}
	\({L^2a^2\over \pi^2}+ (x-y)^2\right)^{-2} .
}$$
Recall that $L=\log T$ and $\tilde\gamma=L\gamma/2\pi$.
And for the sum on $F_2$ we have,
$$\eqalign{
\sum_{0<\gamma,\gamma^\prime,\gamma^{\prime\prime}<T}
	F_2(\gamma,\gamma^\prime,\gamma^{\prime\prime})  =&
\sum_{0<\gamma,\gamma^\prime,\gamma^{\prime\prime}<T}
f_2\(\tilde{\gamma},\tilde{\gamma}^\prime,\tilde{\gamma}^{\prime\prime}\)\cr
=&\mathop{{\sum}'}_{0<\gamma,\gamma^\prime,\gamma^{\prime\prime}<T}
f_2\(\tilde{\gamma},\tilde{\gamma}^\prime,\tilde{\gamma}^{\prime\prime}\)
	+ {2\over 4a^2} \mathop{{\sum}'}_{0<\gamma,\gamma^\prime<T}
f_2\(\tilde{\gamma},\tilde{\gamma}^\prime\) \cr
	&+  \mathop{{\sum}'}_{0<\gamma,\gamma^\prime<T}
f_1\(\tilde{\gamma},\tilde{\gamma}^\prime\)
	+{1\over 2^4a^4} \sum_{0<\gamma<T} 1
}$$
where
$$\eqalign{
f_2(x,y,z)=&{L^4\over 2^4\pi^4}
	\({L^2a^2\over \pi^2}+ (x-y)^2\)^{-1}
	 \({L^2a^2\over \pi^2}+ (y-z)^2\)^{-1}\cr
f_2(x,y)=&{L^2\over 2^2\pi^2}
	 \({L^2a^2\over \pi^2}+ (x-y)^2\right)^{-1} .
}$$
It is interesting that $\sum f_2\(\tilde{\gamma},\tilde{\gamma}^\prime\)$ is 
the sum encountered in [GG].

Now apply $(2.0)$ to the expressions above.  The resulting
integrals can be evaluated explicitly.
With the help of a computer-algebra package we find,
$$\eqalign{
\mathop{{\sum}'}_{0<\gamma,\gamma^\prime,\gamma^{\prime\prime}<T}
f_1\(\tilde{\gamma},\tilde{\gamma}^\prime,\tilde{\gamma}^{\prime\prime}\)
&\sim {3\over 2\pi}T L  \Biggl(
 {{{L^2}}\over {576\,{a^4}}}
- {L\over {128\,{a^5}}} 
+{{23}\over {1536\,{a^6}}} 
- {3\over {256\,{a^7}\,L}} \cr
&\phantom{xxxxxxxxxxxxxx}+T^{-2a}\( 
{1\over {96\,{a^6}\,}} 
+ {1\over {96\,{a^7}\,L}} 
\)\cr
&\phantom{xxxxxxxxxxxxxx}+ T^{-4a}\(
{1\over {1536\,{a^6}}} 
+ {1\over {768\,{a^7}\,L}} 
      \)
\Biggr)
	\cr
\mathop{{\sum}'}_{0<\gamma,\gamma^\prime<T}
f_1\(\tilde{\gamma},\tilde{\gamma}^\prime\)
&\sim {1\over \pi}T L  \Biggl(
 {L\over {64\,{a^3}}}
-{{1}\over {32\,{a^4}}} 
+ {3\over {128\,{a^5}\,L}} 
-T^{-2a}\( {1\over {64\,{a^4}}} 
+ {3\over {128\,{a^5}\,L}} 
\)
 \Biggr)
	\cr
\mathop{{\sum}'}_{0<\gamma,\gamma^\prime,\gamma^{\prime\prime}<T}
f_2\(\tilde{\gamma},\tilde{\gamma}^\prime,\tilde{\gamma}^{\prime\prime}\)
&\sim {3\over 2\pi}T L  \Biggl(
 {{{L^2}}\over {48\,{a^2}}}
- {{5\,L}\over {96\,{a^3}}} 
+{{25}\over {384\,{a^4}}} 
- {7\over {192\,{a^5}\,L}} 
+ T^{-2a}{1\over {24\,{a^5}\,L}} \cr
&\phantom{xxxxxxxxxxxxxx} - T^{-4a}\({1\over {384\,{a^4}}} 
+ {1\over {192\,{a^5}\,L}} 
\)
\Biggr)
	\cr
\mathop{{\sum}'}_{0<\gamma,\gamma^\prime<T}
f_2\(\tilde{\gamma},\tilde{\gamma}^\prime\)
&\sim {1\over \pi}T L  \Biggl(
 {L\over {8\,a}}
-{{1}\over {8\,{a^2}}} 
+ {1\over {16\,{a^3}\,L}} 
- T^{-2a}{1\over {16\,{a^3}\,L}} 
\Biggr)
}$$

Inserting these into our previous formulas gives
$$
\intl_0^T \Sigma_a(t)^3\,dt \sim
{1\over 8 a^3}TL^3\(
 {1} 
+{{3}\over {2\,{a^2}\,L^2}} \)
,
$$
as claimed.  This completes the proof of Proposition~3, and so that of 
Theorem~1.

\no {\bf 5.  Limitations of the method}

The GUE Hypothesis refers to summations over zeros with distinct
indices.  In this paper we encountered unrestricted sums, which
we broke up into several sums over distinct indices.  This step
is not really necessary, for one can write a version of (2.0)
where the summation is not restricted.  In fact, all known methods
of evaluating correlation functions [M] [H] [RS] proceed by first
evaluating an unrestricted sum, and then subtracting out the
``nondistinct'' parts of the sum.  See (4.1) in [RS] for the general
case, and (10) in [H] for a particularly explicit version in the
case of triple correlation.

It is clear that the methods of Section~4 are sufficient,
assuming GUE$N$, to evaluate
$$
\intl_0^T \prod_{n=1}^N \Re\ldz\(\thalf+a_n+it\)\,dt,
\eqno(5.0)
$$
for $a_n\approx 1/\log T$.  The point is that Montgomery's trick works
here, and the rest is just calculating explicit, but messy, integrals.
These last computations could perhaps be simplified by dealing only
with the unrestricted sums described above.
There is nothing intrinsically interesting about the real part
of $\zeta^\prime/\zeta$, but by good fortune the integrand in (1.0)
and in Theorem~1 can be directly related to (5.0) by simple
algebraic identities.  Unfortunately, this trick does not seem to work for
any other integral.  For example, evaluating 
$$
\intl_0^T \lv\ldz\(\thalf+a+it\)\rv^4dt
\eqno(5.1)
$$
by these methods requires evaluation of one of these integrals:
$$
\intl_0^T \({\Im\ldz\(\thalf+a+it\)}\)^4\,dt
\ \ \ \ \  
\ \ \ \ \  
\hbox{or}
\ \ \ \ \  
\ \ \ \ \  
\intl_0^T \({\Re\ldz\(\thalf+a+it\)}\)^2\({\Im\ldz\(\thalf+a+it\)}\)^2\,dt .
$$
The presence of the imaginary part of $\zeta^\prime/\zeta$
makes these much more difficult than (5.0).  In a few cases this 
difficulty can be avoided by simple tricks.  For example, on GUE2,
$$
\intl_0^T\ldz\(\thalf+a+it\)^2\,dt\ll \log^3 T
\ \ \ \ \ \
\ \ \ \ \ \
\hbox{and}
\ \ \ \ \ \
\ \ \ \ \ \
\intl_0^T \(\Re\ldz{\(\thalf+a+it\)}\)^2dt \sim T\,{1-T^{-2a}\over 8a^2} ,
$$
which leads to
$$
\intl_0^T \(\Im\ldz{\(\thalf+a+it\)}\)^2dt \sim T\,{1-T^{-2a}\over 8a^2} .
\eqno(5.2)
$$
Proving (5.2) directly seems to be difficult.  One can use 
the partial fraction formula to give an explicit expression 
for $\Im(\zeta^\prime/\zeta)$, but the methods of Section~4 fail
to yield (5.2).  In particular, Montgomery's trick does not work here.
If these difficulties could be overcome then integrals like (3.1) and (5.1)
could be evaluated directly.  This could lead to a proof that 
(3.0) implies GUE3 for the $\zeta$-function.
It would also give a consistency check between (3.0) and GUE4.
It is interesting that (3.0) can be used to conjecture a value
for (5.1), and this can be used to deduce the value of the two integrals
displayed below  (5.1).

\no {\bf 6. References}

\vskip 6pt

{\parskip 6pt

[BF]  {\sl E.~Bombieri} and {\sl J.B.~Friedlander}, Dirichlet
polynomial approximations to zeta functions, preprint.

[CG]  {\sl J.B.~Conrey} and {\sl A.~Ghosh}, Mean values of the
Riemann zeta-function, III, Proceedings of the Amalfi Conference
On Analytic Number Theory.

[F1] {\sl D.W.~Farmer}, Long mollifiers of the Riemann zeta-function, 
Mathematika {\bf 40} (1993), 88-101.

[F2]  {\sl D.W.~Farmer}, Mean value of Dirichlet series associated with
holomorphic cusp forms, Journal of Number Theory, to appear.

[GM]  {\sl P.X.~Gallagher} and {\sl J.~Mueller}, Primes and zeros
in short intervals, J.~Reine
angew.~Math {\bf 303} (1978), 205-220.

[G1]  {\sl D.A.~Goldston}, On the pair correlation conjecture for
zeros of the Riemann zeta-function, J.~Reine angew.~Math {\bf 385}
(1988), 24-40.

[G2]  {\sl D.A.~Goldston}, On the function $S(T)$ in the theory of
the Riemann zeta-function, Journal of Number Theory {\bf 27} (1987),
149-177.

[GG] {\sl D.A.~Goldston} and {\sl S.M.~Gonek}, Mean values of 
$\zeta^\prime\over \zeta$ and primes in short intervals, preprint

[H] {\sl D.A.~{H}{e}jhal}, On the triple correlation of zeros of
the zeta function, International Mathematics Research Notices {\bf 7}
(1994), 293-302.

[M]  {\sl H.L.~Montgomery}, The pair-correlation of zeros of the
zeta function, Proc.~Symp.~Pure Math 24, 181-193, AMS, Providence 1973.

[Odl]  {\sl A.M.~Odlyzko}, The $10^{20}$ zero of the Riemann zeta function
and 70 million of its neighbors, AT\&T Preprint, 1989.

[RS] {\sl Z.~Rudnick} and {\sl P.~Sarnak}, Zeros of principal
$L$-functions and random matrix theory, preprint.

[T]  {\sl E.C.~Titchmarsh} and {\sl R.~Heath-Brown}, The Theory of the 
Riemann Zeta-Function, Second Ed., Oxford 1986.
}

\vskip 12pt

\vbox{\parskip 0pt
\no Mathematical Sciences Research Institute

\no 1000 Centennial Drive

\no Berkeley, CA \ \ 94720

\no farmer@msri.org
}

\bye